\documentstyle[12pt,twoside,a4mlnew,math30,etamak,etahead,ready]{article}

\def\myleftmark{Lesch and Tolksdorf: On the determinant}
\newcommand{\comment}[1]{\relax}

\def\endproof{\hfill\qed\par\vspace{1em}\bigbreak\noindent}

\begin{document}
\pagestyle{empty}
\renewcommand{\tilde}{\widetilde}
\def\VMark{\relax}

\def\G4-1.3{{\rm (2.3)}}

\renewcommand{\thefootnote}{\arabic{footnote}}
\begin{center}
\LARGE\bf On the Determinant of One-Dimensional Elliptic Boundary Value 
Problems

\bigskip
\bigskip
\Large Matthias Lesch and J\"urgen Tolksdorf

\bigskip\normalsize
Institut f\"ur Mathematik,
Humboldt--Universit\"at zu Berlin,
Unter den Linden 6,
D--10099 Berlin,
Germany,\\
e-mail: lesch@mathematik.hu-berlin.de, tolkdorf@mathematik.hu-berlin.de,\\
Internet: http://spectrum.mathematik.hu-berlin.de/${\sim}$lesch

\bigskip
\bigskip
\today
\end{center}

\begin{abstract}
We discuss the $\zeta-$regularized determinant of elliptic boundary value
problems on a line segment. Our framework is applicable for separated and 
non-separated boundary conditions. 
\end{abstract}

\bigskip\bigskip
\tableofcontents
\newpage\hphantom{++}\newpage
\pagestyle{headings}
\setcounter{page}{1}
\section{Introduction}\label{Sec0}

In \cite{BFK1, BFK2, BFK3}, Burghelea, Friedlander and
Kappeler cal\-cu\-la\-ted the $\zeta-$re\-gu\-la\-rized 
deter\-minant of elliptic differential operators on a line 
segment with periodic and separated boundary 
con\-ditions. In \cite{BFK2} they also discussed 
pseudo\-differential operators over $S^1$, e.g. on $[0,1]$
with periodic boundary conditions. In \cite{L2}, the first 
named author of this paper considered the 
$\zeta-$regularized determinant of second order 
Sturm-Liouville operators with 
regular singularities at the boundary. The common 
phenomenon of \cite{BFK1, BFK2, BFK3} and of
\cite{L2} is that the $\zeta-$regularized determinant is
expressed in terms of a determinant (in the sense of linear 
algebra) of an endomorphism of a finite-dimensional vector 
space of solutions of the corresponding homogeneous 
differential equation.

In this paper we want to show that this phenomenon remains
valid for arbitrary (e.g. non-separated, non-periodic)
boundary conditions and that there exists a simple proof
which works for all types of boundary conditions 
simultaneously. However, our result is less explicit
than the results of \cite{BFK1}, \cite{BFK3} for periodic
and separated boundary conditions, respectively. The 
reason is that for arbitrary boundary conditions we could not 
prove a general deformation result for the variation of the 
leading coefficient. On the other hand while \cite{BFK3}
is limited to even order operators we deal with operators
of arbitrary order (see also the discussion at the end 
of Section 3.1).

The main feature of our approach is the new proof of the
variation formula, Proposition 3.1 below, which uses the
explicit formulas for the resolvent kernel. This together
with some general considerations about $\zeta-$regularized
determinants and regularized limits (Section 2.3) easily give
the main results, Theorem 3.2 and Theorem 3.3.

Since this paper may be viewed as the second part of \cite{L2},
we refer the reader to the end of Section 1 of loc. cit. for a
more detailed historical discussion of $\zeta-$regularized
determinants for one-dimensional operators. Nevertheless, we try 
to keep this paper notationally as self-contained as possible.

The paper is organized as follows: In Section 2
we introduce some notation and review the basic facts 
about the $\zeta-$regularized determinant of an 
operator. In Section 3 we state and prove our main 
results.

This work was supported by Deutsche 
Forschungsgemeinschaft.
   
\section{Generalities}\label{Sec1}

\subsection{Regularized integrals}
First we briefly recall regularized limits and integrals
(c.f. \cite[(1.8)-(1.13c)]{L1}). Let $f:\R\rightarrow\C$ be 
a function having an asymptotic expansion
\begin{equation}
f(x)\sim_{x\rightarrow 0+}\;
\sum_{\Re\alpha\leq 0}x^\alpha\, P_\alpha({\rm log}x) +
o(1)\mylabel{G4-1.1},
\end{equation}
where $P_\alpha\in\C[t]$ are polynomials and $P_\alpha=0$
for all but finitely many $\alpha$. Then we put
\begin{equation}
\LIM_{x\rightarrow 0} f(x) := P_0(0).\mylabel{G4-1.2}
\end{equation}
If $f$ has an expansion like \myref{G4-1.1} as 
$x\rightarrow\infty$ then ${\rm LIM}_{x\rightarrow\infty}$ 
is defined likewise.

Next let $f:\R\rightarrow\C$ such that
\begin{eqnarray}
f(x) &=&
\sum_\alpha x^\alpha\;P_\alpha({\rm log}x) + f_1(x),
\nonumber\\
&=&
\sum_\beta x^\beta\;Q_\beta({\rm log}x) + f_2(x)
\mylabel{G4-1.3},
\end{eqnarray}  
with $P_\alpha, Q_\beta\in\C[t], 
f_1\in{\rm L}^1_{loc}([0,\infty)), 
f_2\in{\rm L}^1([1,\infty))$. Then, we define
\begin{equation}
\regint^\infty_0 f(x)\, dx :=
\LIM_ {a\rightarrow 0}\,\int\limits^1_a f(x)\,dx +
\LIM_ {b\rightarrow\infty}\,\int\limits^b_1 f(x)\,dx
\mylabel{G4-1.4}.
\end{equation}
We note that for a slightly more restricted class of 
functions, this regularized integral can also be defined by 
the Mellin transform (\cite{BS}, \cite[Sec. 2.1]{L1}, 
\cite[(1.12)]{L2}). Note that
\begin{equation}
\regint^\infty_0 x^\alpha\,{\rm log}^k x \, dx = 0
\mylabel{G5-1.5}
\end{equation}
for $\alpha\in\C, k\in\Z_+$ (cf. \cite[(1.13 a-c)]{L2}).

\subsection{Boundary value problems on a line segment}
We consider a linear differential operator
\begin{equation}
l := \sum^n_{k=0} a_k(x)\,D^k, \qquad D:=-i\frac{d}{dx}
\mylabel{G5-1.6}
\end{equation}
defined on the bounded in\-ter\-val $I:=[a,b]$ with ma\-trix 
coefficients $a_k\in C^\infty(I, M(m,\C))$. We assume 
\myref{G5-1.6} to be elliptic, 
i.e. $\det a_n(x)\not=0, x\in I$. A priori, the
differential operator $l$ acts on $C^\infty(I, \C^m)$. We
consider the following boundary condition:
\begin{equation}
{\cal B}(f) := 
R_a\pmatrix{f(a)\cr f'(a)\cr\vdots\cr f^{(n-1)}(a)}
+ 
R_b\pmatrix{f(b)\cr f'(b)\cr\vdots\cr f^{(n-1)}(b)}
= 0\mylabel{G5-1.7},
\end{equation}
where $R_a, R_b\in M(nm, \C)$ are matrices of size 
$nm\times nm$.

We denote by $L:=l_{{\cal B}}$ the differential operator 
\myref{G5-1.6}, restricted to the domain
\begin{equation}
{\cal D}(L) := 
\{f\in H^n(I, \C^m)\;|\; {\cal B}(f)=0\}.\mylabel{G5-1.6a}
\end{equation}
Let $\phi: I\rightarrow M(nm, \C)$ be the fundamental 
matrix of $l$, which means that $\phi$ is the solution of 
the initial value problem
\begin{eqnarray}
\phi'(x) 
&=&
 A\,\phi(x),\nonumber\\[.5em]\mylabel{G6-1.8a}
\phi(a) 
&=& 
{\bf 1}_{nm}\mylabel{G6-1.8b}. 
\end{eqnarray}
Here, $A\in C^\infty(I, M(nm,\C) )$ is the matrix
\begin{equation}
A:=
\pmatrix{0 & {\bf 1}_m & 0 & \ldots & 0\cr
                0 & 0 & {\bf 1}_m & \ldots & 0\cr
        \vdots &\vdots& \ddots &      &\vdots\cr
                0 & 0 & 0 & \ldots & {\bf 1}_m\cr
\beta_0&\beta_1&\beta_2&\ldots&
\beta_{n-1}}\mylabel{G6-1.9},
\end{equation}
where, respectively, 
$\beta_k\equiv -\alpha^{-1}_n\,\alpha_k:=
-(-i)^{-n}\,a^{-1}_n\,(-i)^k\,a_k,\; k=0,\cdots,n-1$, i.e.
$\alpha_k:=(-i)^k\,a_k,\;k=0,\cdots,n$ and 
${\bf 1}_m\in M(m,\C)$ denotes the $m\times m$ unit-matrix.

Sometimes we also write $\phi(x; l)$ to make the dependence 
on the operator \myref{G5-1.6} explicit. We introduce the
matrices
\begin{eqnarray}
{\cal R} &:=&
{\cal R}(l, R_a, R_b) := R_a + R_b\,\phi(b; l),
\mylabel{G6-1.10a}\\[.5em]
{\cal R}(z) &:=& 
{\cal R}(l+z, R_a, R_b) := R_a + R_b\,\phi(b; l+z)
\mylabel{G6-1.10b}.
\end{eqnarray}

It is a well-known fact that the operator $L$ is invertible
if and only if the matrix $\cal R$ is invertible. In this
case the inverse operator $L^{-1}$ is a trace class operator
with kernel
\begin{equation}\mylabel{G6-1.11}
K(x,y) = 
\left\{ \begin{array}{r@{\quad\hbox{if}\quad}l}
-\left[\phi(x)\,{\cal R}^{-1}\,R_b\,\phi(b)\,
\phi^{-1}(y)\right]_{1n}\,\alpha_n(y)^{-1}
& y > x, \\[.8em]
-\left[\phi(x)\,({\cal R}^{-1}\,R_b\,\phi(b) - {\bf 1})\,
\phi^{-1}(y)\right]_{1n}\,\alpha_n(y)^{-1}
& y < x. 
\end{array} \right.
\end{equation}
Here, $[\qquad]_{1n}$ means the upper right entry of a
$n\times n$ block matrix. 
Note that $K(x,y)\in{\bf M}(m,\C)$.

\subsection{The $\zeta-$regularized determinant} 
We briefly discuss $\zeta-$regularized determinants in an
abstract setting. Let ${\cal H}$ be a Hilbert space and let 
$L$ be an (unbounded) operator in ${\cal H}$.

For $\alpha<\beta$ we denote by
\begin{equation}
C_{\alpha,\beta} := 
\{z\in\C\backslash\{0\}\;|\;\alpha\leq\arg z\leq\beta\}
\end{equation}
a sector in the complex plane. We assume that the operator
$L$ has $\theta$ as a principal angle. By this we mean that
there exists an $\epsilon>0$ such that
\begin{equation}
\spec L\cap C_{\theta-\epsilon,\theta+\epsilon}=
\emptyset\mylabel{G7-1.12}.
\end{equation}

Furthermore, we assume that
\begin{equation}
|\!|(L-z)^{-1}|\!|_{{\cal L}({\cal H})}\leq c\,|z|^{-1}, \quad
z\in C_{\theta-\epsilon,\theta+\epsilon},\quad 
|z|\geq{\rm R}\mylabel{G7-1.13}, 
\end{equation}
where $(L-z)^{-1}$ is trace class and there is an asymptotic
expansion in $C_{\theta-\epsilon,\theta+\epsilon}$ as
$z\rightarrow\infty$
\begin{equation}
{\rm Tr}(L-z)^{-1} \sim_{z\rightarrow\infty}
\sum_{\Re\alpha\geq -1-\delta} 
z^\alpha\,P_\alpha(\log z) + o(|z|^{-1-\delta})
\mylabel{G7-1.14},
\end{equation}
where, again, $P_\alpha\in\C[t]$ are polynomials and
$P_\alpha\neq0$ for at most finitely many $\alpha$.

Moreover, we assume that
\begin{equation} 
\deg P_{-1}=0\mylabel{G8-1.15},
\end{equation}
i.e., there are no terms like $z^{-1}\,\log^k(z), k\geq 1$.

The trace class property of $(L-z)^{-1}$ implies that
\begin{equation}
\lim_{{z\rightarrow\infty\atop 
z\in C_{\theta-\delta,\theta+\delta}}}
{\rm Tr}(L-z)^{-1} = 0
\end{equation}
for any $\delta<\epsilon$. Thus, $P_\alpha=0$ if 
$\Re\alpha\geq0$.

In view of \myref{G7-1.13} we can construct the complex
powers of the operator $L$ as follows 
(cf. \cite[Sec. 1]{Se2}, \cite[Sec. 10.1]{Sh}): let 
$\Gamma=\Gamma_1\cup\Gamma_2\cup\Gamma_3$ be the
contour in $\C$ with
\begin{eqnarray}
\Gamma_1 &:=& 
\{r\,e^{i(\theta+2\pi)}\;|\;\rho<r<\infty\},\nonumber\\[.5em]
\Gamma_2 &:=& 
\{\rho\,e^{i(\theta+\varphi)}\;|\;0<\varphi<2\pi\},
\nonumber\\[.5em]
\Gamma_3 &:=& 
\{r\,e^{i\theta}\;|\;\rho<r<\infty\}\mylabel{G8-1.16}.
\end{eqnarray}
Here, the contour $\Gamma$ is traversed such that the set
$\C\backslash\{r\,e^{i\theta}\;|\; r>\rho\}$ lies "inside"
$\Gamma$. Moreover, $\rho$ is chosen so small that
${\rm spec}L\cap\{z\in\C\;|\;|z|\leq\rho\}\subset\{0\}$.

Then, put for $\Re z<0$
\begin{equation}
L_z := \frac{i}{2\pi}\,
\int_\Gamma \lambda^z\,(L-\lambda)^{-1}\,d\lambda
\mylabel{G8-1.17}.
\end{equation}
Here, the complex powers $\lambda^z$ are defined by
$(r\,e^{i(\theta+\varphi)})^z:=r^z\,e^{iz(\theta+\varphi)},
0\leq\varphi\leq 2\pi$.
The same proof as in \cite[Thm.1]{Se1} (cf. also
\cite[Prop. 10.1]{Sh}) shows that $z\mapsto L_z$ is a 
holomorphic semigroup of bounded operators in the Hilbert 
space ${\cal H}$. 

For $k\in\Z, k<0$ we have $L_k=(L_{-1})^k$. Moreover, if
$0\not\in\spec L$ then $L_{-1}=L^{-1}$. If
$0\in\spec L$, then $LL_{-1}$ is a projection onto a
complementary subspace of $\ker L$. Therefore, we
shall write $L^z$ instead of $L_z$.

By assumption $(L-z)^{-1}$ is trace class and in view of
\myref{G7-1.13} we can estimate the trace norm
\begin{eqnarray}
|\!|(L-z)^{-1}|\!|_{\rm tr}
 &\leq&
|\!|(L-z_0)^{-1}|\!|_{\rm tr}\;|\!|(L-z_0)(L-z)^{-1}|\!|
\nonumber\\[.5em]
&\leq&
|\!|(L-z_0)^{-1}|\!|_{\rm tr}\,(1+|z-z_0|\,|\!|(L-z)^{-1}|\!|)
\nonumber\\[.5em]
&\leq&
C,\quad |z|\geq{\rm R}\mylabel{G9-1.18}.
\end{eqnarray}

Therefore, if $\Re z<-1$ the integral \myref{G8-1.17}
converges in the trace norm and the $\zeta-$function
of $L$ with respect to the principal angle $\theta$
\begin{eqnarray}
\zeta_{L,\theta}(s) 
&:=& 
{\rm Tr}(L^{-s}) = 
\sum_{\lambda\in{\rm spec}(L)\backslash\{0\}}
\lambda^{-s}
\nonumber\\[.5em]
&=&
\frac{i}{2\pi}\,
\int_\Gamma z^{-s}\,{\rm Tr}(L-z)^{-1}\,dz
\mylabel{G9-1.19}
\end{eqnarray}
is a holomorphic function for $\Re s>1$.

Furthermore, the asymptotic expansion \myref{G7-1.14}
implies that $\zeta_{L,\theta}(s)$ has a meromorphic
continuation to $\Re s>-\delta$ with poles in the
set $\{\alpha+1\;|\; P_\alpha\neq 0\}$. The order of
the pole $\alpha+1$ is either $\deg P_\alpha$ if 
$\alpha+1\in\Z$ or ${\rm deg}P_\alpha+1$ if 
$\alpha+1\not\in\Z$ (see for instance \cite[Lemma 2.1]{BL}).
Because of the assumption \myref{G8-1.15}
$\zeta_{L,\theta}(s)$ is regular at $0$.

Following Ray and Singer, \cite{RS}, we put 
$\det_\theta L=0$ if $0\in\spec L$, and otherwise
\begin{equation}
{\rm det}_\theta L := \exp(-\zeta'_{L,\theta}(0))
\mylabel{G9-1.20}. 
\end{equation}

It is convenient to deal with the principal angle
$\theta=\pi$. We therefore consider the operator
\begin{equation}
{\tilde L} := e^{i(\pi-\theta)}\,L\mylabel{G10-1.21}.
\end{equation}

Obviously, this operator has $\theta=\pi$ as a principal
angle and it satisfies \myref{G7-1.13}-\myref{G8-1.15},
too. Furthermore,
\begin{equation}
{\tilde L}^{-s} = e^{is(\theta-\pi)}\,L^{-s}\mylabel{G10-1.22},
\end{equation}
and thus
\begin{equation}
\zeta_{{\tilde L},\pi}(s) = 
e^{is(\theta-\pi)}\,\zeta_{L,\theta}(s)\mylabel{G10-1.23}.
\end{equation}

Consequently,
\begin{eqnarray}
\zeta_{L,\theta}(0) 
&=&
\zeta_{{\tilde L},\pi}(0),\nonumber\\[.5em]
\zeta'_{L,\theta}(0)
&=&
\zeta'_{{\tilde L},\pi}(0)+i(\pi-\theta)\,
\zeta_{{\tilde L},\pi}(0)\mylabel{G10-1.24}
\end{eqnarray}
and therefore
\begin{equation}
{\rm det}_\theta L = 
e^{i(\theta-\pi)\,\zeta_{{\tilde L},\pi}(0)}\;
{\rm det}_\pi {\tilde L}\mylabel{G10-1.25}.
\end{equation}

In the sequel we thus assume $\theta=\pi$. We then write
the expansion \myref{G7-1.14} in the form
\def\schnoerkel{{'}}
\addtocounter{equation}{-13}
\begin{equation}
{\rm Tr}(L+x)^{-1} \sim_{x\rightarrow\infty}
\sum_{\Re\alpha\geq -1-\delta} 
x^\alpha\,P_\alpha(\log x) + o(x^{-1-\delta}),\quad
x\geq 0.\mylabel{G7-1.14'}
\end{equation}
\resetschnoerkel\addtocounter{equation}{12}
 
Of course, there exist formulas relating 
the $P_\alpha$ in \myref{G7-1.14} and the corresponding
$P_\alpha$ in \myref{G7-1.14'}.  

\begin{lemma}{S10-1.1}
Let the operator $L$ be given as above with principal angle
$\theta=\pi$. Then,
\begin{eqnarray}
\zeta_{L,\pi}(s) 
&=&
\frac{\sin \pi s}{\pi}\,
\regint^\infty_0 x^{-s}\,
{\rm Tr}(L+x)^{-1}\,dx\mylabel{G10-1.26a},\\[.5em]
\zeta'_{L,\pi}(0) 
&=&
\regint^\infty_0 {\rm Tr}(L+x)^{-1}\,dx\,.
\mylabel{G10-1.26b}
\end{eqnarray}
\end{lemma}

\proof
With respect to the decomposition
${\cal H}=(\ker L)\oplus (\ker L)^\bot$, the operator $L$
reads
\begin{equation}
L = \pmatrix{0 & T \cr
                     0 & L_1},
\end{equation}
where $L_1$ is invertible. In view of \myref{G5-1.5} we
have 
\begin{equation}
\regint^\infty_0 x^{-s}{\rm Tr}(L+x)^{-1}\,dx =
\regint^\infty_0 x^{-s}{\rm Tr}(L_1+x)^{-1}\,dx
\end{equation}
and thus we may assume $L$ to be invertible.

From the estimate \myref{G9-1.18} we conclude that the
following integral is absolutely convergent for
$1<{\rm Re}s<2$:
\begin{eqnarray}
\int_0^\infty x^{-s}\,
[{\rm Tr}(L+x)^{-1}-{\rm Tr}(L^{-1})]\,dx &=&
\sum_{\lambda\in\spec (L)\backslash\{0\}}
\int_0^\infty x^{-s}\,
[(\lambda+x)^{-1}-\lambda^{-1}]\,dx\nonumber\\[.5em]
&=&
\sum_{\lambda\in{\rm spec}(L)\backslash\{0\}}
\regint^\infty_0 x^{-s}\,(\lambda+x)^{-1}\,dx
\nonumber\\[.5em]
&=&
\frac{\pi}{\sin \pi s}\,
\sum_{\lambda\in{\rm spec}(L)\backslash\{0\}}
\lambda^{-s}.
\mylabel{1-2.34}
\end{eqnarray}
Here, we have used \myref{G5-1.5} again.  Hence, the 
first formula is proved.

Since \myref{G7-1.14'}, \myref{G8-1.15} and \cite[(1.12)]{L1}
imply
\begin{equation}
\regint^\infty_0 x^{-s}\,{\rm Tr}(L+x)^{-1}\,dx =
\frac{P_{-1}(0)}{s} +
\regint^\infty_0 {\rm Tr}(L+x)^{-1}\,dx +
O(s),\quad s\rightarrow 0
\mylabel{1-2.35}
\end{equation}
we reach the conclusion by noting that 
$\frac{\sin \pi s}{\pi}=s+O(s^3), s\rightarrow 0$.
\endproof

\begin{lemma}{S12-1.2}
Let $L$ be as before, $\theta=\pi$. Then, we have the
asymptotic expansion
\begin{equation}
\log{\rm det}_\pi (L+x) \sim_{x\rightarrow\infty}
\sum_{\Re\alpha\geq-1-\delta}
x^{\alpha+1}\,Q_\alpha(\log x) + O(x^{-\delta})
\mylabel{G12-1.27},
\end{equation}
where $P_\alpha=(\alpha+1)\,Q_\alpha + Q'_\alpha$.
Furthermore, $Q_{-1}(\log x)= P_{-1}(0)\,\log x$. In
particular
\begin{equation}
\LIM_{x\rightarrow\infty}\,\log{\rm det}_\pi (L+x) = 0
\mylabel{G12-1.28}.
\end{equation}
\end{lemma}

\proof
Since $L^{-1}$ is trace class, it follows that
$\log{\rm det}_\pi(L+x)$ is differentiable and
\begin{equation}
\frac{d}{dx}\,\log{\rm det}_\pi (L+x) = 
{\rm Tr}(L+x)^{-1}.
\end{equation}
Hence,
\begin{equation}
\LIM_{y\rightarrow\infty}\,\log{\rm det}_\pi (L+y) -
\log{\rm det}_\pi (L+x) =
\regint^\infty_x {\rm Tr}(L+y)^{-1}\,dy.
\end{equation}
Comparing this equation for $x=0$ with the preceding 
lemma yields \myref{G12-1.28}. Hence,
\begin{eqnarray}
\log{\rm det}_\pi(L+x) 
&=&
-\regint^\infty_x {\rm Tr}(L+y)^{-1}\,dy
\nonumber\\[.5em]
&\sim_{x\rightarrow\infty}&
\sum_{\Re\alpha \geq -1-\delta}
-\regint^\infty_x y^\alpha\,P_\alpha(\log y)\,dy +
O(x^{-\delta})\mylabel{G13-1.28}
\end{eqnarray}
and we reach the conclusion.
\endproof

The reader might ask why we argued so complicated
in order to get the first equality of \myref{G13-1.28}. It
appears to be a direct consequence of \myref{G10-1.26b} via 
the apparently "trivial" calculation
\begin{eqnarray}
\log{\rm det}_\pi(L+x) 
&=&
-\regint^\infty_0 {\rm Tr}(L+x+y)^{-1}\,dy
\nonumber\\[.5em]
&=&
-\regint^\infty_x {\rm Tr}(L+y)^{-1}\,dy.
\end{eqnarray}
However, in general for functions f like \myref{G4-1.3}  
we have
\begin{equation}
\regint^\infty_0 f(x+y)\,dy \neq
\regint^\infty_x f(y)\,dy
\mylabel{G13-1.29}.
\end{equation}
Consequently, some care must be in order. Since
the operator $L^{-1}$ is trace class, the phenomenon
\myref{G13-1.29} does not occur for ${\rm Tr}(L+x)^{-1}$.
More precisely, if $f\in{\rm L}^1_{loc}([0,\infty))$ and
\begin{equation}
f(x)\sim_{x\rightarrow\infty}
\sum_\alpha x^\alpha\,P_\alpha(\log x),
\end{equation}
then
\begin{equation}
\regint^\infty_0 f(x+y)\,dy\,-\,\regint^\infty_x f(y)\,dy =
\LIM_{b\rightarrow\infty}\int^{b+x}_b f(y)\,dy
\end{equation}
and in general this vanishes only if $P_\alpha = 0$ for 
$\alpha\in\Z_+$. As an illustrative example we consider
$f(x):=x^\alpha, \alpha\in\Z$. Then, we get
\begin{equation}
\regint^\infty_0 (x+y)^\alpha\,dy = 
\left\{ \begin{array}{r@{\quad\hbox{if}\quad}l}
-\frac{x^{\alpha+1}}{\alpha+1}
& \alpha\in\Z_-\backslash\{-1\}, \\[.8em]
-\ln x 
& \alpha = -1, \\[.8em]
0 
& \alpha\in\Z_+;
\end{array} \right.
\end{equation}
\begin{equation}
\regint^\infty_x f(y)\,dy =
\left\{ \begin{array}{r@{\quad\hbox{if}\quad}l}
-\frac{x^{\alpha+1}}{\alpha+1}
& \alpha\neq-1, \\[.8em]
-\ln x 
& \alpha = -1.
\end{array} \right.
\end{equation}
Hence,
\begin{equation}
\regint^\infty_0 f(x+y)\,dy\,-\,\regint^\infty_x f(y)\,dy =
\left\{ \begin{array}{r@{\quad\hbox{if}\quad}l}
\frac{x^{\alpha+1}}{\alpha+1}
& \alpha\in\Z_+, \\[.8em]
0 
&\alpha\notin\Z_+ .
\end{array} \right.
\end{equation}

\section{Main results}\label{Sec2}

From now on we restrict ourselves to boundary 
value problems on a line segment as introduced in Sec. 2.2.
Let $(l, {\cal B})$ be an elliptic boundary value problem,
$L:=l_{\cal B}$. More precisely, we assume that $(l, {\cal B})$ 
is elliptic in the sense of \cite[Def.1]{Se1} and that it satisfies
Agmon's condition \cite[Def.2]{Se1}. Agmon's condition
assures that the coefficient $a_n(x)$ has a certain 
principal angle, $\theta$. Then we can find an angle 
$\theta'$, arbitrary close to $\theta$, such that $\theta'$ is a
principal angle for $L$ and $a_n(x)$. Henceforth we shall 
write $\theta$ for a common principal angle of $a_n(x)$ and 
$L$. In short: we will refer to an operator $L=l_{\cal B}$, 
defining an elliptic boundary value problem $(l,{\cal B})$, as 
an {\sl admissible} operator.\\

\subsection{Operators of order $\geq 2$}
If in addition $n\geq 2$ then the conditions 
\myref{G7-1.13}-\myref{G8-1.15} are fulfilled by the work of Seeley 
\cite[Se2]{Se1}. Namely, \myref{G7-1.13} follows from 
\cite[Lemma 15]{Se1} and by \cite[Thm.2]{Se2} we have an asymptotic 
expansion as $z\rightarrow\infty$ in 
$C_{\theta-\epsilon,\theta+\epsilon}$
\begin{equation}
{\rm Tr}(L-z)^{-1}\sim_{z\rightarrow\infty}
\sum_{k=0}^{\infty} a_k\,z^{\frac{1-k}{n}-1}.
\mylabel{V1-1}
\end{equation}
\myref{G8-1.15} is automatically fulfilled since there are no 
log-terms in \myref{V1-1}.\\

Summing up, we see that $\det_\theta L$ is well-defined for $n\geq 2$.
First order operators are slightly more complicated since in
this case $(L-z)^{-1}$ is not of trace class. This problem
will be treated separately in subsection 3.2.

First, we study the behavior of $\det_\theta L$ under
deformations of the coefficients of $l$.

\begin{prop}{S16-2.1}
Assume that the coefficients $a_0,\ldots, a_{n-2}$ depend
smoothly on a parameter $t$. Let $L_t$ be the
corresponding family of operators. If $L_t$ is invertible
then we have
\begin{equation}
\partial_t\,{\rm log}{\rm \det}_\theta\,L_t =
\partial_t\,{\rm log}{\rm det}{\cal R}_t\mylabel{V1-2}
\end{equation}
where ${\cal R}_t:={\cal R}(l_t, R_a, R_b)$, 
cf. \myref{G6-1.10a}.
\end{prop}
\proof  The inclusion 
$H^2([a,b], \C^m)\hookrightarrow L^2([a,b], \C^m)$ is trace
class and ${\cal D}(L)\subset H^n([a,b], \C^m)$.
Hence the operators $D^k L^{-1}$ are trace class, as well,
for $k=0,\ldots,n-2$.  Hence,
\begin{eqnarray}
\partial_t\,{\rm log}{\rm det}_\theta L_t &=&
{\rm Tr}((\partial_t L_t)L^{-1}_t)\nonumber\\[.5em]
&=&
\sum_{j=0}^{n-2}\int^b_a{\rm tr}_{\C^m}
\left(\partial_t a_j(t;x)\,(D^j K_t)(x,x)\right)\,dx
\nonumber\\[.5em]
&=&
\sum_{j=0}^{n-2}\int^b_a{\rm tr}_{\C^m}
\left(\partial_t \alpha_j(t;x)\, K^{(j)}_t(x,x)\right)\,dx
\nonumber\\[.5em]
&=&
\sum_{j=0}^{n-2}\int^b_a{\rm tr}_{\C^m}
\left(-\alpha^{-1}_n(x)\partial_t \alpha_j(t;x)\, 
[{\tilde K}^{(j)}_t(x,x)]_{1n}\right)\,dx
\nonumber\\[.5em]
&=&
\sum_{j=0}^{n-2}\int^b_a{\rm tr}_{\C^m}
\left(\partial_t \beta_j(t;x)\,
[{\tilde K}^{(j)}_t(x,x)]_{1n}\right)\,dx,
\end{eqnarray}
where
\begin{equation}\mylabel{V1-30} 
{\tilde K}_t(x,y) := 
\left\{ \begin{array}{r@{\quad\hbox{if}\quad}l}
\phi_t(x)\,{\cal R}_t^{-1}\,R_b\,\phi_t(b)\,
\phi_t^{-1}(y)
& y > x, \\[.8em]
\phi_t(x)\,({\cal R}_t^{-1}\,R_b\,\phi_t(b) - {\bf 1})\,
\phi_t^{-1}(y)
& y < x. 
\end{array} \right.
\end{equation}
Note that ${\tilde K}^{(j)}_t$ is continuous on the diagonal
for $j=0,\ldots,n-2$, but ${\tilde K}^{(n-1)}_t$ has a jump.
This is one of the reasons that this proposition is limited
to the case of constant $a_{n-1}$.
Thus we have
\begin{eqnarray}
\partial_t\,{\rm log}{\rm det}_\theta L_t &=&
{\rm Tr}((\partial_t L_t)\,L^{-1}_t)
\nonumber\\[.5em]
&=&
\int^b_a{\rm tr}_{\C^{nm}}
\left((\partial_t A_t)(x)\,{\tilde K}_t(x,x)\right)\,dx.
\end{eqnarray}

By use of \myref{G6-1.8a} we then calculate:
\begin{eqnarray}
{\rm tr}_{\C^{nm}}
\left((\partial_t A_t)(x)\,{\tilde K}_t(x,x)\right) &=&
{\rm tr}_{\C^{nm}}\left[
(\partial_t A_t)(x)\,\phi_t(x)\,{\cal R}_t^{-1}\,
R_b\,\phi_t(b)\,\phi_t^{-1}(x)\right]
\nonumber\\[.5em]
&=&
{\rm tr}_{\C^{nm}}\left[
(\partial_x\partial_t \phi_t)(x)\,{\cal R}_t^{-1}\,
R_b\,\phi_t(b)\,\phi_t^{-1}(x)\right]
\nonumber\\[.5em]
& &
-\left[(\partial_x\phi_t)(x)\,\phi^{-1}_t(x)\,
(\partial_t \phi_t)(x)\,{\cal R}_t^{-1}\,
R_b\,\phi_t(b)\,\phi_t^{-1}(x)\right]
\nonumber\\[.5em]
&=&
\partial_x{\rm tr}_{\C^{nm}}\left[
(\partial_t \phi_t)(x)\,{\cal R}_t^{-1}\,
R_b\,\phi_t(b)\,\phi_t^{-1}(x)\right]
\mylabel{V1-50} 
\end{eqnarray}
to obtain
\begin{eqnarray}
\partial_t\,{\rm log}{\rm det}_\theta L_t &=&
{\rm tr}_{\C^{nm}}\left[(\partial_t\phi_t)(b)
\,{\cal R}_t^{-1}\, R_b\right]
\nonumber\\[.5em] 
&=&
\partial_t\,{\rm log}{\rm det}{\cal R}_t,
\end{eqnarray}
which proves the statement.\endproof

\begin{theorem}{S16-2.2}
Let $(l, {\cal B})$ be an admissible operator of order $n\geq 2$,
$L:=l_{{\cal B}}$. Assume the principal angle $\theta$ equals
$\pi$. For ${\cal R}(z):={\cal R}(l+z, {\cal B})$ we obtain an
asymptotic expansion
\begin{equation}
{\rm log}{\rm det}{\cal R}(z)\sim_{z\rightarrow\infty}
\sum_{{k=0\atop k\neq 1}}^\infty
b_ k\,z^{\frac{1-k}{n}} + b_1 + \zeta_{L,\pi}(0)\,{\rm log}z
\end{equation}
in a conic neighborhood of $\R_+$. Furthermore,
\begin{equation}
{\rm log}{\rm det}_\pi\,(L+z) =
{\rm log}{\rm det}\,{\cal R}(z) -
\LIM_{w\rightarrow\infty}\,
{\rm log}{\rm det}\,{\cal R}(w).
\end{equation}
\end{theorem}

\proof
In view of \myref{V1-1} we have an asymptotic expansion
\begin{equation}
{\rm Tr}(L+z)^{-1}\sim_{z\rightarrow\infty}
\sum_{k=0}^{\infty} a_k\,z^{\frac{1-k}{n}-1}.
\end{equation}
We apply the preceding proposition with
$a_0(z; x)=a_0(x)+z$ and $a_k(z; x)=a_k(x), k\geq 1$. Then
\begin{eqnarray}
\partial_z\,{\rm log}{\rm det}{\cal R}(z) 
&=&
\partial_z\,{\rm log}{\rm det}_\pi\,(L + z)
\nonumber\\[.5em]
&=&
{\rm Tr}(L + z)^{-1}
\nonumber\\[.5em]
&\sim_{z\rightarrow\infty}&
\sum_{k=0}^\infty a_k\,z^{\frac{1-k}{n}-1}.
\mylabel{V1-3}
\end{eqnarray}
This proves the first assertion. Note that from Lemma 2.1
one easily concludes $a_1=\zeta_{L,\pi}(0)$.  By \myref{V1-3} 
${\rm log}{\rm det}{\cal R}(z)-{\rm log}{\rm det}_\pi(L+z)$
is a constant. Then the second assertion follows from
\myref{G12-1.28}.
\endproof

\begin{theorem}{V18-2.2}
Let again $(l, {\cal B})$ be an admissible operator of order
$n\geq 2$,
$L:=l_{{\cal B}}$, with principal angle $\theta=\pi$. We put
\begin{equation}
\log{\cal C}(l,{\cal B}):= 
-\LIM_{z\rightarrow\infty}\log{\rm det}{\cal R}(z).
\mylabel{V1-0}
\end{equation}
Then,
\begin{equation}
{\rm det}_\pi L = {\cal C}(l,{\cal B})\det{\cal R},\quad
{\cal R}:={\cal R}(0).
\mylabel{V1-0'}
\end{equation}
Furthermore, ${\cal C}(l,{\cal B})$ depends only on 
$a_n, a_{n-1}$ and the boundary operator ${\cal B}$, i.e. 
${\cal C}(l,{\cal B})={\cal C}_1(a_n,a_{n-1},R_a,R_b)$.
\end{theorem}

\proof
\myref{V1-0} and \myref{V1-0'} are immediate
consequences of the preceding Theorem. To prove the
last statement we consider two admissible operators
\begin{equation}
l_j:=\sum_{k=0}^n a_{k,j}(x)\,D^k,\quad j=0,1\,,
\nonumber
\end{equation}
where $a_{n,0}=a_{n,1}, a_{n-1,0}=a_{n-1,1}$
and the boundary condition ${\cal B}$ is fixed.

We put 
$$l_t:=t\,l_1 + (1-t)\,l_0,\quad 0\leq t\leq 1.$$
We would like to apply Proposition 3.1. However, it may 
happen that 
${\rm spec}\,l_t\cap\{z\in\C\,|\,z\leq 0\}\neq\emptyset$ 
for some $t$. But since $\pi$ is a principal angle
for the leading symbol of $l_t$ there exists a $z_0>0$
such that $L_t+z$ is invertible for all $z\geq z_0$.

By Proposition 3.1 we then have
${\cal C}(l_0+z, {\cal B})={\cal C}(l_1+z,{\cal B})$ for 
$z>z_0$. Since both functions are holomorphic we are done.
\endproof

Note that formulas \myref{V1-0} and \myref{V1-0'} express
the $\zeta-$regularized determinant of $L$ completely in
terms of the solutions of the homogeneous differential
equation $(L+z)u=0$. It seems impossible, however, to find
an explicit formula for the coefficient ${\cal C}(l,{\cal B})$ 
in full generality. But in cases where the fundamental 
matrix ${\cal R}(z)$ can be calculated explicitly one can also
find an expression for ${\cal C}(l,{\cal B})$.

Now we are going to discuss in detail non-separated 
boundary conditions for second order operators. We
therefore consider the following\\

\noindent
{\bf Example}:
Let $A, B, C, D\in M(m,\C)$ and consider the operator
$$l := -\frac{d^2}{dx^2} + q(x)$$
with boundary operator ${\cal B}=(R_a,R_b)$, where
$$R_a := \pmatrix{A & B  \cr
                C & D  \cr},\quad R_b := {\bf 1}_{2m}.$$
It turns out that the operator $L=l_{\cal B}$ is admissible iff 
the meromorphic function
\begin{equation}
M(z):=
\det(-z\frac{B}{2}+\frac{A+D}{2}-\frac{1}{z}\,\frac{C}{2})
\mylabel{V1-4}
\end{equation}
does not vanish identically. Hence, let us assume $M(z)\neq0$.
Note that $M$ is a Laurent polynomial.

By the preceding proposition 
${\cal C}(l, {\cal B})=:{\cal C}_2(A,B,C,D)$ is independent of 
$q$. More precisely,
\begin{prop}{V19-2.3}
${\cal C}(l, {\cal B})^{-1}$ is equal to the leading coefficient 
of the Laurent polynomial $M(z)$; $\zeta_{L,\pi}(0)$ equals
$\frac{1}{2}$ the degree of $M(z)$.
\end{prop}

\proof
As remarked before it suffices to consider the case $q=0$. 
Then the fundamental solution $\phi(x,z):=\phi(x, l+z)$ reads
\begin{equation}
\phi(x,z)=
\pmatrix{\cosh[(x-a)\sqrt{z}]\,{\bf 1}_m &  
\frac{\displaystyle\sinh[(x-a)\sqrt{z}]}
{\displaystyle\sqrt{z}}\,{\bf 1}_m\cr
\sqrt{z}\,\sinh[(x-a)\sqrt{z}]\,{\bf 1}_m & 
\cosh[(x-a)\sqrt{z}]\,{\bf 1}_m},\nonumber
\end{equation}
where, again, ${\bf 1}_m$ denotes the $m\times m$ unit 
matrix.

For the rest of the proof all matrices will be $2\times 2$
block matrices with $m\times m$ block entries and for
simplicity we will omit ${\bf 1}_m$. Abbreviating
$c:=b-a, w:=\sqrt{z}$ we find
\begin{eqnarray}
\phi(b,z) &=&
\pmatrix{\cosh cw & 
\frac{\displaystyle\sinh cw}{\displaystyle w}\cr
w\,\sinh cw & \cosh cw},\nonumber\\[.5em]
&=&
W\,\pmatrix{e^{cw} &    0 \cr
                          0    & e^{-cw}}\,W^{-1},\nonumber
\end{eqnarray}
where
$$W = \pmatrix{1 & 1 \cr
                          w & -w}.$$
Thus,
\begin{eqnarray}
\det{\cal R}(z) &=& \det(R_a+R_b\phi(b,z))\\[.5em]
&=&
\det(R_a+W\,\pmatrix{e^{cw} &    0 \cr
                          0    & e^{-cw}}\,W^{-1})\\[.5em]
&=&
\det(W^{-1\,}R_a\,W+\pmatrix{e^{cw} &    0 \cr
                          0    & e^{-cw}})\\[.5em]
&=&
e^{mcw}\,{\rm det}_{\C^m}[W^{-1}\,R_a\,W]_{22}
+ O(w^{m+1}\,e^{(m-1)cw}),\nonumber
\end{eqnarray}
since 
$$W^{-1}\,R_a\,W = w\,X_1 + X_0 + w^{-1}\,X_{-1},
\quad X_i\in M(2m,\C),\,\, i=-1,0,1.$$
Here, 
$$[W^{-1}\,R_a\,W]_{22} = -\frac{w}{2}\,B
+ \frac{A+D}{2} - \frac{1}{2w}\,C$$
denotes the lower right entry of the $2\times 2$ block 
matrix $W^{-1}\,R_a\,W$.

This implies
\begin{eqnarray}
\log{\rm det}{\cal R}(z) &=&
mc\sqrt{z} + 
\log\!\!\left[{\rm det}_{\C^m}[W^{-1}\,R_a\,W]_{22}
+ O(z^{\frac{m+1}{2}}\,e^{-c\sqrt{z}})\,\right]
\nonumber\\[.5em]
&=&
mc\sqrt{z} + \log M(\sqrt{z}) + 
O\left(\hbox{$\frac{z^{\frac{m+1}{2}}\,
e^{-c\sqrt{z}}}{M(\sqrt{z})}$}\right).\nonumber
\end{eqnarray}
Since $M(w)$ is a Laurent polynomial we may write
$$M(\sqrt{z}) = \lambda\,z^{k/2}+O(z^{\frac{k-1}{2}}),
\quad z\rightarrow\infty$$
and thus
$$\log M(\sqrt{z}) = \frac{k}{2}\,\log z + 
\log\lambda + O(z^{-\frac{1}{2}})
\quad z\rightarrow\infty$$
and we reach the conclusion.
\endproof

The leading coefficient of $M(z)$ is in general difficult to
describe. Thus, it seems hard to find a more explicit formula
for ${\cal C}(l,{\cal B})$ than given in the preceding
Proposition.

We discuss some special cases:
\begin{enumerate}
\item $B$ invertible:
\begin{eqnarray}
M(z) &=& 
(-1)^m\,(\det\frac{B}{2})\,z^m\,\det({\bf 1}_m
+O(z^{-1}))\nonumber\\[.5em]
&=&
(-1)^m\,(\det\frac{B}{2})\,z^m + O(z^{m-1}));\nonumber
\end{eqnarray}
\item $B=0,\, A+D$ invertible:
$$M(z) = \det(\frac{A+D}{2}) + O(z^{-1});$$
\item $B=0, A+D=0$:
$$M(z) = \frac{(-1)^m}{2^m}\,\det C\,z^{-m}.$$
\end{enumerate}
Hence,
\begin{equation}
{\cal C}_2(A,B,C,D) =
\left\{ \begin{array}{r@{\quad\hbox{if}\quad}l}
\frac{\displaystyle (-1)^m\,2^m}{\displaystyle \det B}
& \det B\neq 0, \\[.9em]
\frac{\displaystyle 2^m}{\displaystyle \det(A+D)}
& B=0,\, \det(A+D)\neq 0,\\[.9em]
\frac{\displaystyle (-1)^m\,2^m}{\displaystyle \det C}
& B=0,\, A+D=0. 
\end{array} \right.\mylabel{V1-5}
\end{equation}
Of course, this does not cover all possible cases. The periodic
boundary conditions are given by $R_a=-{\bf 1}_{2m}$
and thus
$${\cal C}_2(A,B,C,D)=(-1)^m,$$
which is consistent with \cite[Thm.1]{BFK1}.\\

Next, we discuss how ${\cal C}(a_n, a_{n-1}, {\cal B})$ 
depends on the coefficients $a_n$ and $a_{n-1}$. We start 
with the dependence on the subleading coefficient $a_{n-1}$
and use the standard trick to eliminate it (cf. also
\cite[Prop.2.2]{BFK3}).

For this let again $L=l_{\cal B}$ be an admissible operator, 
${\cal B}=(R_a, R_b)$. Let also $U:I\rightarrow M(m,\C)$ be 
the unique solution of the initial value problem
\begin{eqnarray}
U'(x) &=& -\frac{i}{n}\,(a_n^{-1}a_{n-1})(x)\,U(x)
\nonumber\\[.5em]
U(a) &=& {\bf 1}_m\,.\mylabel{V1-6}
\end{eqnarray}
The determinant of $U(x)$ is given by
\begin{equation}
\det U(x) =
\exp\left(-\frac{i}{n}\,\int^x_a
{\rm tr}(a_n^{-1}\,a_{n-1})(y)\,dy\right)\neq0.
\end{equation}
By conjungation of $l$ with $U$ we find
\begin{equation}
l^u:=U^{-1}\,l\,U = \sum_{j=0}^n {\tilde a}_j(x)\,D^j
\end{equation}
with ${\tilde a}_n(x)=(U^{-1} a_n U)(x)$ and 
${\tilde a}_{n-1}=0$. Since 
${\rm spec}\,{\tilde a}_n = {\rm spec}\,a_n$, ${\tilde a}_n$ 
has the same principal angle as $a_n$. Furthermore, for
$L^u:=U^{-1}LU$ we have
\begin{equation}
\spec L^u = \spec L
\end{equation}
and hence $\det_\theta L^u=\det_\theta L$.
Next, we determine the transformed boundary operator
${\cal B}^u:=(R_a^u,R_b^u)$. If 
\begin{equation}
\phi = 
\pmatrix{\varphi_1&\cdots&\varphi_n\cr
                \vdots &\phantom{\cdots}&\vdots\cr
                \varphi_1^{(n-1)}&\cdots&\varphi_n^{(n-1)}}
\end{equation}
denotes a fundamental matrix of $l$, then the corresponding
fundamental matrix of $l^u$ reads
\begin{eqnarray}
{\tilde\phi}^u &=& 
\pmatrix{U^{-1}\varphi_1&\cdots&U^{-1}\varphi_n\cr
\vdots &\phantom{\cdots}&\vdots\cr
(U^{-1}\varphi_1)^{(n-1)}&\cdots&(U^{-1}\varphi_n)^{(n-1)}}
\nonumber\\[.8em]
&=&
T(U)\,\phi,
\end{eqnarray}
where
\begin{equation}
T(U)_{ij}(x) := 
\left\{ \begin{array}{r@{\quad\hbox{if}\quad}l}
{i\choose j}\,(\partial_x^{(i-j)}\,U^{-1})(x)
& 0\leq j\leq i\leq n-1, \\[.5em]
0\quad\quad
& j>i\,.
\end{array} \right.\mylabel{V1-6'}
\end{equation}
Note that $\det T(U)=(\det U)^{-n}$. 

However, the fundamental matrix ${\tilde\phi}^u$ is not 
normalized. We therefore put
\begin{eqnarray}
\phi^u(x) &:=& 
{\tilde\phi}^u(x)\,({\tilde\phi}^u)^{-1}(a)\nonumber\\[.5em]
&=&
T(U)(x)\,\phi(x)\,T(U)^{-1}(a).
\end{eqnarray}

We now determine the boundary conditions for $L^u$. Let
$g\in{\cal D}(L^u)$. Then, $g=U^{-1}f$ with 
$f\in{\cal D}(L)$. With $F:=(f,\ldots f^{(n-1)})^t,
G:=(g,\ldots g^{(n-1)})^t$ we have
\begin{equation}
G=T(U)\,F
\end{equation}
and thus
\begin{eqnarray}
0 &=& R_a\,F(a) + R_b\,F(b)\nonumber\\[.5em]
&=&
R_a\,T(U)^{-1}(a)\,G(a)+R_b\,T(U)^{-1}(b)\,G(b)
\nonumber\\[.5em]
&=:&
{\tilde R}^u_a\,G(a) + {\tilde R}^u_b\,G(b).
\end{eqnarray}

We put
\begin{eqnarray}
R^u_a &:=& T(U)(b)\,{\tilde R}^u_a\nonumber\\[.5em]
R^u_b &:=& T(U)(b)\,{\tilde R}^u_b.
\end{eqnarray}
Note that $R^u_a, R^u_b$ define the same boundary 
conditions as ${\tilde R}^u_a, {\tilde R}^u_b$.             
 
Then,
\begin{eqnarray}
{\cal R}(l^u+z, R^u_a, R^u_b) &=&
R^u_a + R^u_b\,\phi^u(b, l^u+z)\nonumber\\[.5em]
&=&
T(U)(b)\,{\cal R}(l+z, R_a, R_b)\,T(U)^{-1}(a)
\end{eqnarray}
and thus
\begin{equation}
\det{\cal R}(l^u+z, R^u_a, R^u_b) =
(\det U(b))^{-n}\,\det{\cal R}(l+z, R_a, R_b).
\end{equation}

Consequently,
\begin{equation}
{\rm det}_\theta L^u = {\rm det}_\theta L 
\end{equation}
implies
\begin{equation}
{\cal C}(a_n, a_{n-1}, R_a, R_b) =
(\det U(b))^{-n}\,{\cal C}(U^{-1}a_nU, 0, R^u_a, R^u_b).
\mylabel{V1-7}
\end{equation}

We thus have proved the 
\begin{prop}{V19-2.4} Let 
$L=l_{\cal B}$ be an admissible operator and let $U(x)$ be 
the unique solution of the initial value problem \myref{V1-6}. 
Then, the operator $L^u=l^u_{{\cal B}^u}$ is also admissible. 
It has the same principal angle $\theta$ as $L$ and 
\begin{eqnarray}
{\cal C}(a_n, a_{n-1}, R_a, R_b) =
\exp\left(i\,\int^b_a{\rm tr}(a_n^{-1}\,a_{n-1})(y)\,dy\right)
\,{\cal C}(U^{-1}a_nU, 0, R^u_a, R^u_b),
\end{eqnarray}
with
\begin{eqnarray}
R^u_a &:=& T(U)(b)\,R_a\,T(U)^{-1}(a),\\[.5em]
R^u_b &:=& T(U)(b)\,R_b\,T(U)^{-1}(b)
\mylabel{V1-8}
\end{eqnarray}
and $T(U)$ is defined by formula \myref{V1-6'}.
\end{prop}

As an application, we consider the operator
$$l := -\frac{d^2}{dx^2} + p(x)\,\frac{d}{dx} + q(x)$$
with the same boundary operator ${\cal B}=(R_a, R_b)$ as
given in the preceding example. Again, by 
Proposition 3.3 it is sufficent to consider $q=0$. Notice that
the leading coefficient $a_2$ of the operator $l$ is invariant 
with respect to conjungation with $U$. Hence, in the two 
specific cases where, respectively, $B$ is invertible, or 
$B=0$ and $A+D$ is invertible, we can simply make use of 
\myref{V1-5} to obtain 
\begin{equation}
{\cal C}_2(p, A, B, C, D) =
e^{\hbox{$\frac{1}{2}$}\int^b_a {\rm tr}\,p(x)\,dx}
\left\{ \begin{array}{r@{\,\,\hbox{if}\,\,}l}
\frac{ (-1)^m\,2^m}{ \det B}
& \det B\neq 0, \\[.8em]
\frac{ 2^m}{\det(A+D)}
& B=0,\, \det(A+D)\neq 0. 
\end{array} \right.
\mylabel{V1-9}
\end{equation}\\

We now turn to the dependence of  
${\cal C}(l, {\cal B})$ on the leading symbol of the 
differential operator $l$.  The aim is to get an 
explicit formula analogous to \myref{V1-7} in the case of 
$a_{n-1}$. Unfortunately, this is much more involved
than in the case of separated boundary conditions. 
Following \cite{BFK3} one considers the family
of operators $l_t:=\alpha_t\,(D^n+l')$, where $l'$ 
denotes a differential operator of order $n-1$ 
and $\alpha_t(x),\,t\in [0,1]$ is a smooth variation
of $a_n(x)$ such that $\alpha_0 =\hbox{Id}$ and
$\alpha_1=a_n$. Then, the question arises whether
the corresponding operators $L_t:=(l_t, {\cal B})$ 
are admissible for all $t\in[0,1]$. To answer this
question seems to be hopeless for the general situation 
discussed in this paper. Note, however, for
a given admissible operator $L=(l,{\cal B})$ the constant 
${\cal C}(l, {\cal B})$ can be calculated if the fundamental 
solution of the corresponding homogeneous equation is 
known. 

\subsection{Operators of order 1}
In this subsection we briefly indicate how Theorem 3.2 and 3.3
generalize to operators of order one. Let $L=l_{\cal B}$ be an
admissible operator of order one with principal angle $\pi$. A priori
$(L+x)^{-1}$ is not of trace class. However, the trace of
$(L+x)^{-1}$ can be regularized.

In the sequel we use the notation ${\rm Res}_k\,f(z_0)$ for the
coefficient of $(z-z_0)^k$ in the Laurent expansion of the 
meromorphic function $f$.

The function ${\rm Tr}(L+z)^{-s}$ is meromorphic with simple
poles in $1,0,-1,\ldots$, which follows from \myref{3.1-4}
below, and we put
\begin{equation}
{\overline{\Tr}}(L+z)^{-1}:={\rm Res}_0\,(L+z)^{-s}|_{s=1}.
\mylabel{3.1-1} 
\end{equation}
Then,
\begin{eqnarray}
\frac{d}{dz}{\rm log}{\det}_{\pi}\,(L+z) 
&=&
-\frac{d}{ds}|_{s=0}\,\frac{d}{dz}\,{\rm Tr}(L+z)^{-s}
\nonumber\\[.5em]
&=&
\frac{d}{ds}|_{s=0}\,s{\rm Tr}(L+z)^{-s-1}
\nonumber\\[.5em]
&=&
{\overline{\Tr}}(L+z)^{-1}
\nonumber\\[.5em]
&=&
{\rm Tr}(\,(L+z)^{-1}-L^{-1}\,)+{\overline{\Tr}}(L^{-1}),
\mylabel{3.1-2}
\end{eqnarray}
since $(L+z)^{-1}-L^{-1} = -z(L+z)^{-1}L^{-1}$ is of trace class.

The same calculation as \myref{1-2.34} shows that
\begin{eqnarray}
\zeta_{L,\pi}(s) 
&=&
\frac{\sin \pi s}{\pi}\,
\regint^\infty_0 x^{-s}\,{\rm Tr}[(L+x)^{-1}- L^{-1}]\,dx
\nonumber\\[.5em]
&=&
\frac{\sin \pi(s-1)}{\pi(s-1)}\,
\regint^\infty_0 x^{1-s}\,{\rm Tr}(L+x)^{-2}\,dx.
\mylabel{3.1-3}
\end{eqnarray}

By the work of Seeley we have an asymptotic expansion 
\begin{equation}
{\rm Tr}(L+x)^{-2}\sim_{x\rightarrow\infty}
\sum_{k=0}^\infty a_k\,x^{-k-1},
\mylabel{3.1-4}
\end{equation}
which  implies in view of $\sin\pi(s-1)=\pi(s-1)^2+ O((s-1)^3)$ that
\begin{equation}
{\overline{\Tr}}(L^{-1})= {\rm Res}_0\,\zeta_{L,\pi}(1) =
\regint^\infty_0 {\rm Tr}(L+x)^{-2}\,dx.
\mylabel{3.1-5}
\end{equation}

On the other hand, we have
\begin{eqnarray}
{\overline{\Tr}}(L^{-1})=\regint^\infty_0 {\rm Tr}(L+x)^{-2}\,dx
&=&
\LIM_{R\rightarrow\infty}\left\{-\int_0^R
\frac{d}{dx}{\rm Tr}[(L+x)^{-1}-L^{-1}]\,dx\right\}
\nonumber\\[.5em]
&=&
-\LIM_{R\rightarrow\infty}{\rm Tr}[(L+R)^{-1}-L^{-1}].
\mylabel{3.1-6}
\end{eqnarray}

Comparing this with \myref{3.1-2} gives
\begin{equation}
\LIM_{z\rightarrow\infty}\frac{d}{dz}{\rm log}{\rm det}_{\pi}(L+z)
=\LIM_{z\rightarrow\infty}{\overline{\Tr}}(L+z)^{-1}=0.
\mylabel{3.1-7}
\end{equation}

Summing up we can state the anolog of Lemma 2.1 and 2.2 for operators of
first order:

\begin{lemma}{3.1-l3.6}
Let $L$ be as before. Then we have
\begin{eqnarray}
\zeta_{L,\pi}(s) 
&=&
\frac{\sin \pi s}{\pi}\,
\regint^\infty_0 x^{-s}\,{\overline{\Tr}}(L+x)^{-1}\,dx
\mylabel{3.1-8},\\[.5em]
\zeta'_{L,\pi}(0) 
&=&
\regint^\infty_0 {\overline{\Tr}}(L+x)^{-1}\,dx\,.
\mylabel{3.1-9}
\end{eqnarray}
and the asymptotic expansions in a conic neighborhood of $\R_+$
\begin{equation}
{\overline{\Tr}}(L+x)^{-1}
\sim_{x\rightarrow\infty}
\sum^\infty_{k=1}\frac{a_k}{k}\,x^{-k} + a_0\log x,
\mylabel{3.1-10}
\end{equation}
\begin{equation}
\log{\rm det}_\pi(L+x)
\sim_{x\rightarrow\infty}
\sum_{k=2}^{\infty}\frac{a_k}{k(k-1)}\,x^{1-k} + \zeta_{L,\pi}(0)\,\log x
+ a_0\,x\log x + a_0\,x,
\mylabel{3.1-11}
\end{equation}
in particular
\begin{eqnarray}
\LIM_{x\rightarrow\infty}\,{\overline{\Tr}}(L+x)^{-1} &=& 0,
\mylabel{3.1-12}\\[.5em]
\LIM_{x\rightarrow\infty}\,\log{\rm det}_\pi\,(L+x) &=& 0.
\mylabel{3.1-13} 
\end{eqnarray}
\end{lemma}

\proof
The equation \myref{3.1-8} follows from \myref{3.1-3} and \myref{3.1-9}
follows from \myref{3.1-8}, similar to \myref{1-2.35}; \myref{3.1-10},
\myref{3.1-11} follow from integrating the expansion \myref{3.1-4};
\myref{3.1-12} follows from \myref{3.1-7} and finally, \myref{3.1-13}
is proved exactly as \myref{G12-1.28}.
\endproof

\begin{theorem}{3.1-t3.7}
Let $L$ be as before. For ${\cal R}(x)={\cal R}(l+x,{\cal B})$ we have
an asymptotic expansion
\begin{equation}
\log\det{\cal R}(x)\sim_{x\rightarrow\infty}
\sum_{k=2}^{\infty}\frac{a_k}{k(k-1)}\,x^{1-k} + \zeta_{L,\pi}(0)\,\log x
+ b +  a_0\,x\log x + c\,x.
\mylabel{3.1-14}
\end{equation}
Furthermore,
\begin{equation}
\log{\rm det}_\pi\,(L+x)-\log\det{\cal R}(x)= -b-(c-a_0)x.
\mylabel{3.1-15}
\end{equation}
\end{theorem}

\proof
We cannot apply Proposition \plref{S16-2.1} to $L+z$ since the operator $L$ is of
order one. But the same computation as in the proof of Proposition 
\plref{S16-2.1}
shows that
\begin{eqnarray}
\frac{d}{dx}{\overline{\Tr}}(L+x)^{-1} 
&=&
-{\rm Tr}(L+x)^{-2}
\nonumber\\[.5em]
&=&
-\int^b_a{\rm tr}(L+x)^{-2}(t,t)\,dt
\nonumber\\[.5em]
&=&
\frac{d}{dx}\,\int^b_a{\rm tr}[(L+x)^{-1}(t,t+0)]\,dt
\nonumber\\[.5em]
&=&
\frac{d^2}{dx^2}\,\log\det{\cal R}(x),
\mylabel{3.1-16}
\end{eqnarray}
hence,
\begin{equation}
\log{\rm det}_\pi\,(L+x)-\log\det{\cal R}(x)
\mylabel{3.1-17}
\end{equation}
is a polynomial of degree one and we are done.
\endproof

As a consequence, we end up with the formula
\begin{equation}
{\rm det}_\pi\,L = e^{-b}\,\det{\cal R},
\mylabel{3.1-18}
\end{equation}
where $b:=\LIM_{x\rightarrow\infty}\log\det{\cal R}(x)$.

\begin{small}
\def\and{{\rm and }}

\end{small}


\begin{thebibliography}{GrSe12}


\bibitem[BFK1]{BFK1}
{\sc D. Burghelea, L. Friedlander and T. Kappeler:}
\newblock {\it On the determinant of elliptic differential
and finite difference operators in vector bundles over $S^1$.}
\newblock Commun. Math. Phys. {\bf 138} (1991), 1-18

\bibitem[BFK2]{BFK2}
{\sc D. Burghelea, L. Friedlander and T. Kappeler:}
\newblock {\it Regularized determinants for 
pseudodifferential operators in vector bundles over $S^1$.}
\newblock Int. Equ. Op. Th. {\bf 16} (1993), 496-513

\bibitem[BFK3]{BFK3}
{\sc D. Burghelea, L. Friedlander and T. Kappeler:}
\newblock {\it On the determinant of elliptic boundary
value problems on a line segment.}
\newblock Proc. Am. Math. Soc. {\bf 123}, No. 10, (1995)
3027-3038 

\bibitem[BL]{BL}
{\sc J. Br\"uning, M. Lesch:}
\newblock {\it On the spectral geometry of algebraic curves.}
\newblock J. reine angew. Math. {\bf 474} (1996), 25-66

\bibitem[BS]{BS}
{\sc J. Br\"uning, R. Seeley:}
\newblock {\it Regular Singular Asymptotics.}
\newblock Adv. Math. {\bf 58} (1985), 133-148

\bibitem[L1]{L1}
{\sc M. Lesch:}
\newblock {\it Operators of Fuchs type, conical singularities
and asymptotic methods.}
\newblock Teubner Texte zur Mathematik, vol. 136,
 B. G. Teubner, Leipzig, 1997

\bibitem[L2]{L2}
{\sc M. Lesch:}
\newblock {\it Determinants of regular singular 
Sturm-Liouville operators.}
\newblock To appear in Math. Nachr.

\bibitem[RS]{RS}
{\sc D. B. Ray, I. M. Singer:}
\newblock {\it R-torsion and the Laplacian on Riemannian
manifolds.}
\newblock Adv. Math. {\bf 7} (1971), 145-210

\bibitem[Se1]{Se1}
{\sc R. Seeley:}
\newblock {\it The resolvent of an elliptic boundary problem.}
\newblock Amer. J. Math. {\bf 91} (1969), 889-920

\bibitem[Se2]{Se2}
{\sc R. Seeley:}
\newblock {\it Analytic extension of the trace associated
with elliptic boundary problems.}
\newblock Amer. J. Math. {\bf 91} (1969), 963-983 

\bibitem[Sh]{Sh}
{\sc M. A. Shubin:}
\newblock {\it Pseudodifferential operators and spectral
theory.}
\newblock Springer, Berlin-Heidelberg-New York, 1987


\end{thebibliography}
\end{document}